\documentclass[letterpaper, 10 pt, conference]{ieeeconf}  

\IEEEoverridecommandlockouts                              
\usepackage{hyperref}

\usepackage{amsmath}
\usepackage{amssymb}
\usepackage{accents}
\usepackage{mathtools}
\usepackage{subfig}
\usepackage{comment}
\usepackage{booktabs}
\usepackage{threeparttable}
\usepackage[textwidth=1.5cm, textsize=small]{todonotes}
\usepackage{array}

\newcommand{\ubar}[1]{\underaccent{\bar}{#1}}

\title{\LARGE \bf
Stochastic Model Predictive Control based on Mixed Random Variables for Economic Energy Management
}

\author{
\authorblockN{
Janik Pinter\authorrefmark{1},
Maximilian Beichter\authorrefmark{1},
Ralf Mikut\authorrefmark{1},
Veit Hagenmeyer\authorrefmark{1},
and Frederik Zahn\authorrefmark{1}
}
\authorblockA{
Institute for Automation and Applied Informatics\\
Karlsruhe Institute of Technology (KIT), Germany\\
\authorrefmark{1}\texttt{\{firstname.lastname\}@kit.edu}
}
}

\begin{document}
\bstctlcite{BSTcontrol}

\maketitle
\thispagestyle{empty}
\pagestyle{empty}

\begin{abstract}
Optimal scheduling of batteries has significant potential to reduce electricity costs and to enhance grid resilience. However, effective battery scheduling must account for both physical constraints as well as uncertainties in consumption and generation of renewable energy sources. Instead of optimizing fixed battery power setpoints, we propose an approach that optimizes battery power intervals, allowing the optimization to explicitly account for uncertain consumption and generation as well as how the battery system should respond to them within its physical limits. Our method is based on mixed random variables, represented as mixtures of discrete and continuous probability distributions. Building on this representation, we develop an analytical stochastic formulation for minimizing electricity costs in a residential setting with load, photovoltaics, and battery storage. We demonstrate its effectiveness across real-world data from 15 residential buildings over five consecutive months. Compared with deterministic and probabilistic benchmark controllers, the proposed interval-based optimization achieves the lowest costs. These results show that mixed random variables are a practical and promising tool for decision-making under uncertainty.
\end{abstract}

\section{Introduction}
\label{sec:introduction}
Energy storage technologies such as batteries play an increasingly important role in managing the uncertainty introduced by distributed renewable energy resources. 
In residential buildings with photovoltaics (PV), batteries are often installed to increase self-consumption and to benefit from time-varying electricity prices. 
In this setting, operating batteries efficiently and economically under uncertainty typically relies on optimization-based scheduling, which combines forecasts of consumption, generation and electricity prices in an optimization model to minimize energy costs. 
However, those forecasts contain inherent uncertainty, which has to be considered in the design of the scheduling approach. 
Many state-of-the-art solutions compute optimal schedules as battery power setpoints. 
Applying these setpoints directly, however, neglects forecast errors and the resulting deviations from the expected operating conditions, potentially leading to higher costs.
Furthermore, the optimized schedule must respect physical constraints such as charging power limits and battery capacity.
In the present work, we show that mixed random variables provide a powerful tool for addressing these challenges by enabling the optimization of battery power intervals rather than single setpoints, while accounting for both uncertainty and constraints within a coherent analytical framework.
We introduce a novel formulation based on mixed random variables for cost-optimal battery scheduling under uncertainty. 
This approach integrates scheduling and application into a single optimization problem based on probabilistic forecasts. 
Using real-world data from 15 residential buildings, we show that our approach achieves lower costs than fixed-setpoint and rule-based control approaches.

Battery scheduling under uncertainty is a well-established topic in the energy management literature, and many different approaches have been proposed to model and handle uncertainty. 
In most existing works, the scheduling framework computes single battery power setpoints over the scheduling horizon, as in \cite{lechl_uncertainty-aware_2025,su_optimal_2021}. Since forecast errors are unavoidable in practice, these setpoints are either applied directly while handing the deviation to the grid \cite{bird_lifetime_2025}, or they are complemented by an additional downstream control layer that determines how deviations are handled during operation \cite{pasqui_self-dispatching_2025}. Moreover, uncertain variables appearing in constraints are typically addressed either through chance constraints \cite{lechl_uncertainty-aware_2025} or through deterministic reformulations, for example in robust or interval optimization \cite{dong_hybrid_2023}. Our work follows a different idea:
Instead of scheduling single battery power setpoints, we optimize battery power intervals together with policies that determine the realized battery setpoint during operation as uncertainty unfolds. This embeds the operational response directly into the original optimization problem and removes the need for a separate downstream optimization step.
We note that, although our approach optimizes intervals, it does not correspond to the usual notion of \textit{interval optimization} in the literature \cite{dong_hybrid_2023, su_interval_2020}. There, intervals are typically used to replace uncertain inputs, whereas our framework takes probability density functions (PDFs) of uncertain variables as inputs and yields intervals as outputs. 
We achieve this by employing mixed random variables. In our previous work \cite{pinter_probabilistic_2025}, we introduced the use of mixed random variables in energy applications to derive probabilistic dispatch schedules. In the present work, we extend the underlying theory for electricity cost minimization and embed it in a Model Predictive Control (MPC) framework, shifting the focus from probabilistic day-ahead scheduling to actual cost performance.

The remainder of the paper is organized as follows. 
\autoref{sec:model-derivation} motivates the problem and develops the analytical formulation of the proposed approach, which is then embedded into a stochastic optimization problem in \autoref{sec:optimization}. \autoref{sec:model-simulation} describes the simulation setup used for validation. The results are presented in \autoref{sec:results} and are discussed in \autoref{sec:discussion}. Finally, \autoref{sec:conclusion} concludes the paper.

\section{Model Derivation}
\label{sec:model-derivation}

\begin{figure}[h]
    \centering
    \includegraphics[width=0.52\linewidth]{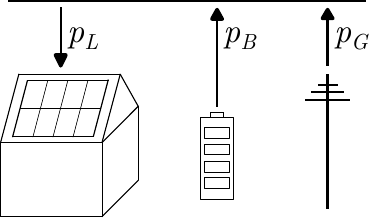}
    \caption{General setting. Net-load $p_L(k)$, battery power $p_B(k)$ and grid power $p_G(k)$ are in balance \cite{pinter_probabilistic_2025}.}
    \label{fig:power-balance}
\end{figure}

We develop our approach for a single residential building equipped with a PV and battery system, as in \cite{su_interval_2020, pinter_probabilistic_2025}. \autoref{fig:power-balance} sketches our setup, where net-load $p_L(k)$, battery power $p_B(k)$ and grid power $p_G(k)$ are in balance:
\begin{equation}
\label{eq:power-balance}
    p_L(k) = p_B(k) + p_G(k).
\end{equation}
In this work, time is discretized into time intervals with length $\Delta t = 1h$, indexed by $k\in\mathcal{K}$, where $\mathcal{K} = \{1, ..., K\}$ denotes the set of discrete time steps. Thus, all powers $p_{i}(k)$ represent power averages over the $k$-th time interval.

\subsection{Problem Motivation}
\label{sec:problem-motivation2}

Many battery control approaches rely on forecasts of the upcoming inflexible net-load to compute battery power schedules. These schedules determine the grid exchange, which is directly responsible for electricity costs. In practice, however, the net-load is inherently uncertain, meaning that predicted and realized net-loads differ. In those instances, the power balance in Equation \eqref {eq:power-balance} must still hold, even though the schedules for $p_B(k)$ and $p_G(k)$ were derived for a different net-load realization. A practical question arises: how should the controller react when the realized net-load differs from the expectation?

One approach is to execute the scheduled battery setpoint and adjust the grid power to satisfy the power balance. As a result, the battery power is known in advance and the grid power becomes uncertain. We refer to this strategy as \textbf{Fixed-Battery}. It is implemented deterministically in \cite{bird_lifetime_2025} and probabilistically in \cite{su_interval_2020}.

Another approach is to adjust the battery power so that the scheduled grid exchange can be maintained \cite{pinter_averaging_2025}. We refer to this strategy as \textbf{Fixed-Grid}. 
Ideally, this makes the grid power known in advance and shifts the uncertainty to the battery power. In practice, however, physical battery limits prevent full compensation of all net-load realizations. As a result, the grid power must compensate for some net-load realizations, making both battery and grid power uncertain.

In the present work, we propose a probabilistic framework that follows the Fixed-Grid idea in a more flexible way: the compensation for deviations is shared between the battery and the grid, based on battery power intervals that are optimized for electricity cost. 
These results cannot be represented by standard probability distributions. Mixed random variables, however, provide exactly this capability by combining probability masses at specific values with continuous probability densities over ranges. 


\subsection{Uncertainty Modeling}

Future net-load $p_L$ is inherently uncertain and can only be estimated.\footnote{For brevity, the time index $k$ is omitted for the remainder of this section.} In this work, we model the net-load $p_L$ as a random variable $P_L:\Omega\rightarrow\mathbb{R}$ defined on a probability space ($\Omega,\mathcal{F}, \mathbb{P}$). $\Omega$ describes the set of real-valued possible outcomes, $\mathcal{F}$ is a sigma algebra of measurable events, and $\mathbb{P}$ is a probability measure. 
$P_L$ is assumed to be absolutely continuous, meaning that a valid Probability Density Function (PDF) $f_{P_L}$ exists.
To denote a realization $\omega \in \Omega$ of a random variable $P_X$, we write $p_X^{\omega}$.

\subsection{Battery Model}
\label{sec:battery-model}

\begin{figure*}[!t]
    \centering
    \subfloat[Net-Load PDF]{%
        \includegraphics[height=3cm]{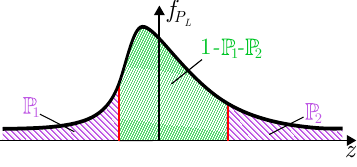}%
        \label{fig:pdf-prosumption}%
    }\hfill
    \subfloat[Battery Power PDF]{%
        \includegraphics[height=3cm]{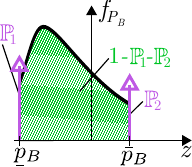}%
        \label{fig:pdf-battery}%
    }\hfill
    \subfloat[Grid Power PDF] {%
        \includegraphics[height=3cm]{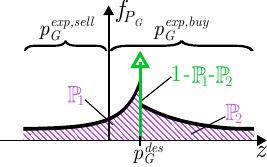}%
        \label{fig:pdf-grid}%
    }
    \caption{Separation of the continuous net-load PDF (a) in two mixed random variable PDFs featuring discrete and continuous properties (b) and (c). The battery power PDF (b) is restricted to a closed interval. The grid power PDF (c) contains a discrete event at $p_G^{des}$. The bigger the interval $[\protect\ubar{p}_B, \Bar{p}_B]$, the more likely it is that the battery can be controlled such that the desired grid exchange $p_G^{des}$ can be followed, thus the bigger the discrete green event in (c). Adapted from \cite{pinter_probabilistic_2025}.}
    \label{fig:pdfs}
\end{figure*}

We derive a probabilistic battery model following the Fixed-Grid strategy. The model adjusts the battery power such that the realized grid exchange is close to a desired grid exchange $p_G^{des}$.
To explicitly reflect this rule, we introduce two parameters $\ubar{p}_B$ and $\Bar{p}_B$. They define the battery power interval that is used to follow $p_G^{des}$ during a specific time step. For any net-load realization $p_L^{\omega}$, perfect tracking of $p_G^{des}$ would require the battery power $p_B^{\omega} = p_L^{\omega} - p_G^{des}$. If this required theoretical value lies outside $[\ubar{p}_B, \Bar{p}_B]$, we deliberately clip the batter action to the nearest bound. 
Accordingly, three operating conditions arise:
\begin{itemize}
    \item If tracking $p_G^{des}$ requires "too much charging" (i.e., a theoretical battery power below $\ubar{p}_B$), we set the battery power to $p_B^{\omega}=\ubar{p}_B$.
    \item If tracking $p_G^{des}$ is feasible within the chosen battery range $p_B^{\omega}\in[\ubar{p}_B, \Bar{p}_B]$, we can set $p_B^{\omega}=p_L^{\omega} - p_G^{des}$ (which yields $p_G^{\omega}=p_G^{des}$ according to Equation \eqref{eq:power-balance}). 
    \item If tracking $p_G^{des}$ requires "too much discharging" (i.e., a theoretical battery power above $\Bar{p}_B$), we set the battery power to $p_B^{\omega}= \Bar{p}_B$.
\end{itemize}
These operating conditions can be formulated probabilistically by introducing a mixed random variable, i.e., the probabilistic battery power $P_B$:
\begin{equation}  
\label{eq:mrv-pb}
    P_{B}(P_L) =
     \begin{cases}
       \ubar{p}_B  &  \ubar{p}_B \ge P_{L} - p_G^{des} \\
        P_{L} - p_G^{des} & \ubar{p}_B < P_{L} - p_G^{des} < \Bar{p}_B  \\
       \Bar{p}_B & \Bar{p}_B \le  P_{L} - p_G^{des}.
     \end{cases} \\
\end{equation}
The PDF\footnote{To be precise, a \textit{PDF} is fully continuous, whereas a \textit{PMF (probability mass function)} is fully discrete. A mixed distribution does not have a commonly used standard term, which is why we refer to the mixed distribution function of a mixed random variable as a PDF for simplicity.} of the mixed random variable then is
\begin{equation}
\label{eq:pdf-pb}
\begin{split}
f_{P_B}(z) ={}& \mathbb{P}_1 \delta(z-\ubar{p}_B) + \mathbb{P}_2 \delta(z-\Bar{p}_B) \\
+ &f_{P_L}(z+p_G^{des})\left[U(z-\ubar{p}_B) - U(z-\Bar{p}_B)\right],
\end{split}
\end{equation}
where $\delta(\cdot)$ is the Dirac delta function and $U(\cdot)$ is the unit step function that restricts the domain in which the continuous part is defined. Additionally, 
\begin{subequations}
\label{eq:p1-p2}
\begin{align}
    \mathbb{P}_1&\coloneqq\mathbb{P}(P_B=\ubar{p}_B) = \mathbb{P}(P_L\le\ubar{p}_B+p_G^{des}) \\
    \mathbb{P}_2&\coloneqq\mathbb{P}(P_B=\Bar{p}_B) = \mathbb{P}(P_L\ge\Bar{p}_B+p_G^{des})
\end{align}
\end{subequations}
describe the probabilities of operating the battery system at its chosen bounds, i.e., $\ubar{p}_B$ or $\Bar{p}_B$ respectively.\footnote{The calculation of the probabilities is done by evaluating the CDF $F_{P_L}$ via $\mathbb{P}_1=F_{P_L}(\ubar{p}_B+p_G^{des})$ and $\mathbb{P}_2=1-F_{P_L}(\Bar{p}_B+p_G^{des})$.} 
See \cite{shankar_probability_2021} or our previous work \cite{pinter_probabilistic_2025} for further information on how to perform mixed transformations of random variables. An exemplary PDF $f_{P_B}(z)$ is sketched in \autoref{fig:pdfs}b. It is neither fully continuous nor fully discrete, which classifies $P_B$ as a \textit{mixed random variable}. 
The PDF shows the likelihood of operating the battery at different power levels (x-axis) such that the desired grid exchange $p_G^{des}$ can be followed.

With Equation \eqref{eq:pdf-pb}, we establish a probabilistic representation of the battery power. To incorporate this mixed random variable into an optimization framework in a tractable way, we approximate it by its expected value, thereby replacing the full uncertainty description with a single representative quantity:
\begin{align}
\label{eq:pb-exp}
        \mathbb{E}[P_B] &= \int_{-\infty}^{\infty} z f_{P_B}(z) dz \nonumber\\
        &=\mathbb{P}_1\ubar{p}_B + \int_{\ubar{p}_B}^{\Bar{p}_B} z f_{P_L} (z+p_G^{des}) dz + \mathbb{P}_2\Bar{p}_B.
\end{align}
With that, we are able to split the expected battery power $\mathbb{E}[P_B]$ into complementary charging and discharging components:\footnote{The complementary constraints in Equations \eqref{eq:pb-exp-relax} and later in \eqref{eq:pg-des-relax} are relaxed to enable nonlinear solvers to handle the problem without introducing binary variables.}
\begin{subequations}
    \begin{align}
    \mathbb{E}[P_B] = p_B^{ch} + p_B^{dis},\quad  &p_B^{ch} \le 0,\quad p_B^{dis} \ge 0, \\
      p_B^{ch} p_B^{dis} &\ge -10^{-8}. \label{eq:pb-exp-relax}
\end{align}
\end{subequations}
This allows us to formulate the battery state evolution according to 
\begin{equation}
    e(k+1) = e(k) - p_B^{ch} \eta^{ch} \Delta t - p_B^{dis}/\eta^{dis} \Delta t.
\end{equation}
$\Delta t$ denotes the time between two steps, $\eta^{ch}$ and $\eta^{dis}$ denote the charging and discharging efficiency respectively, and $e(k)$ denotes the State-of-Energy (SoE) of the battery system.

To summarize, in our approach, the decision variables $\ubar{p}_B$, $\Bar{p}_B$, and $p_G^{des}$ first determine the expected battery power $\mathbb{E}[P_B]$, which is then used to update the SoE $e$. However, enforcing feasibility only for the resulting expected SoE is not sufficient, because $\mathbb{E}[P_B]$ does not necessarily coincide with the battery power that is ultimately realized. Instead, any battery power realization $p_B^{\omega}\in [\ubar{p}_B, \Bar{p}_B]$ may occur and must satisfy the physical battery constraints. Accordingly, we formulate the battery constraints as follows:
\begin{subequations}
\label{eq:physical-bat-constraints}
\begin{align}
    p_B^{\min} &\le \ubar{p}_B \le \Bar{p}_B \le p_B^{\max} \\
    e^{\min} &\le e(k) \le e^{\max} \\
    \ubar{p}_B(k) &\ge (e(k) - e^{\max}) / (\eta^{ch} * \Delta t) \\
    \Bar{p}_B(k) &\le (e(k) - e^{\min}) \eta^{dis} / \Delta t,
\end{align}
\end{subequations}
where $p_B^{\min}$ and $p_B^{\max}$ denote the physical minimum and maximum battery power, while $e^{\min}$ and $e^{\max}$ define the lower and upper bounds of the SoE.

\subsection{Grid Model}
\label{sec:grid-model}

In this section, we derive the probabilistic grid model of the Fixed-Grid strategy.
We assume that a battery system adjusts its power output within $[\ubar{p}_B, \Bar{p}_B]$ to allow the grid to follow a desired value $p_G^{des}$. 
Accordingly, three operating conditions arise. They can be formalized by modeling the probabilistic grid power $P_G$ as a mixed random variable:
\begin{equation}
\label{eq:mrv-pg}
         P_{G}(P_L) = 
     \begin{cases}
       P_{L} - \ubar{p}_B  &\ubar{p}_B \ge  P_{L} - p_G^{des} \\
       p_G^{des} & \ubar{p}_B < P_{L} - p_G^{des} < \Bar{p}_B  \\
       P_{L} - \Bar{p}_B & \Bar{p}_B \le P_{L} - p_G^{des}. \\ 
     \end{cases} 
\end{equation}
Consider for example the upper case of Equation \eqref{eq:mrv-pg}: if tracking $p_G^{des}$ requires "too much charging" (i.e., a theoretical battery power below $\ubar{p}_B$), the battery power is set to its lowest bound $p_B^{\omega} = \ubar{p}_B$ such that the difference between the realized and the desired grid exchange is minimized. In this instance, the grid power follows $p_G^{\omega}=p_L^{\omega}-\ubar{p}_B\neq p_G^{des}$.

With that, the PDF of $P_G$ can be formulated as follows:
\begin{align}
\label{eq:pdf-pg}
    f_{P_G}(z) &= f_{P_L}(z+\ubar{p}_B) [1-U(z-p_G^{des})] \nonumber \\
    &+ (1-\mathbb{P}_1-\mathbb{P}_2)\delta(z-p_G^{des}) \\
    &+ f_{P_L}(z+\Bar{p}_B) U(z-p_G^{des}), \nonumber
\end{align}
with $\mathbb{P}_1$ and $\mathbb{P}_2$ according to Equation \eqref{eq:p1-p2}. The PDF is sketched in \autoref{fig:pdfs}c, where the discrete event $p_G^{\omega}=p_G^{des}$ occurs with a non-zero probability. 
The continuous events (purple) show that deviations from $p_g^{des}$ can still occur with a probability of $\mathbb{P}_1$ and $\mathbb{P}_2$ respectively.
In fact, the shown continuous parts are simply the uncertainty tails of the original net-load uncertainty in \autoref{fig:pdfs}a. 
Thus, the battery system essentially cuts the green part of the net-load uncertainty out, maps the cut-out uncertainty to the discrete event, and shifts the remaining net-load uncertainty tails horizontally. 

Similar to the probabilistic battery power, Equation \eqref{eq:pdf-pg} provides a probabilistic representation of the grid power.
To incorporate the grid power into an optimization framework, we again approximate it by its expected value:
\begin{align}
    \mathbb{E}[P_G] = &\int_{-\infty}^{\infty} z f_{P_G}(z) dz \nonumber\\
                    = &\int_{-\infty}^{p_G^{des}} z f_{P_L}(z+\ubar{p}_B) dz \nonumber \\
                    + &(1-\mathbb{P}_1-\mathbb{P}_2) p_G^{des}  \\
                    + &\int_{p_G^{des}}^{\infty} z f_{P_L}(z+\Bar{p}_B) dz. \nonumber
\end{align}
The expected grid power $\mathbb{E}[P_G]$ is directly linked to the expected electricity costs. 
However, since buying electricity from and selling electricity to the grid are associated with different prices in the residential sector, we further divide the expected grid power into two parts, i.e.,
\begin{subequations}
\label{eq:pg-exp-split}
\begin{align}
    p_G^{{exp, sell}} &= \int_{-\infty}^{p_G^{des, sell}} z f_k(z+\ubar{p}_B) dz\nonumber\\
    &+ (1-\mathbb{P}_{1}-\mathbb{P}_{2}) p_G^{des, sell}  \\ 
    &+ \int_{p_G^{des, sell}}^{0} z f_k(z+\Bar{p}_B) dz \nonumber\\
    p_G^{{exp, buy}} &= \int_{0}^{p_G^{des, buy}} z f_k(z+\ubar{p}_B) dz \nonumber\\
    &+ (1-\mathbb{P}_{1}-\mathbb{P}_{2}) p_G^{des, buy} \\ 
    &+ \int_{p_G^{des, buy}}^{\infty} z f_k(z+\Bar{p}_B) dz, \nonumber
\end{align}
\end{subequations}
where $p_G^{{exp, sell}}$ and $p_G^{{exp, buy}}$ denote the expected exported and imported power. They are not complementary, i.e., both can be nonzero at the same time. Instead, they represent a decomposition of the distribution in \autoref{fig:pdfs}c into two expected values: one associated with the negative part of the distribution and one associated with the positive part. To enable this decomposition, the desired grid exchange $p_G^{des}$ is split into two complementary variables according to:
\begin{subequations}
\label{eq:pg-des-split}
\begin{align}
    &p_G^{des} = p_G^{des, sell} + p_G^{des, buy},  \\
    &p_G^{des, sell} \le 0, \quad p_G^{des, buy} \ge 0,  \\
    &p_G^{des, sell} p_G^{des, buy} \ge -10^{-8}. \label{eq:pg-des-relax}
\end{align}
\end{subequations}
This allows the expected imported and exported energy to be evaluated separately and therefore enables the use of asymmetric import and export prices.

\subsection{Probabilistic Power Balance}

We close this theoretical section by drawing attention to the fundamental power balance of Equation \eqref{eq:power-balance}. We have constructed the probabilistic battery and grid power in Equations \eqref{eq:mrv-pb} and \eqref{eq:mrv-pg} in such a way that our probabilistic power balance holds per definition, i.e.,
\begin{align}
\label{eq:power-balance-prob}
    P_B + P_G &= 
         \begin{cases}
       \ubar{p}_B + P_L - \ubar{p}_B  &  \ubar{p}_B \ge P_{L} - p_G^{des} \\
        P_{L} - p_G^{des} + p_G^{des} & \ubar{p}_B < P_{L} - p_G^{des} < \Bar{p}_B  \\
       \Bar{p}_B + P_L - \Bar{p}_B & \Bar{p}_B \le  P_{L} - p_G^{des}
     \end{cases} \nonumber \\
     &= P_L.
\end{align}
$P_B$ and $P_G$ are dependent random variables since they are linked via $p_G^{des}$, $\ubar{p}_B$ and $\Bar{p}_B$. They are the decision variables of our optimization framework, derived in the following section.

\section{Stochastic Model Predictive Control}
\label{sec:optimization}

We embed the proposed model into a receding-horizon stochastic optimization framework for minimizing residential electricity costs over a planning horizon $\mathcal{K}$. 
The framework takes probabilistic net-load forecasts as input, represented by the PDF $f^k_{P_L}$ and CDF $F^k_{P_L}$ for each time step $k\in\mathcal{K}$. 
The optimization problem is formulated as follows:
\begin{equation}
\begin{aligned}
\min_{\ubar{p}_{B},\, \Bar{p}_{B},\, p_{G}^{des}}\ 
& \sum_{k \in \mathcal{K}} \left(
c^{{buy}}(k)\, p_G^{{exp,buy}}(k)
+ c^{{sell}}(k)\, p_G^{{exp,sell}}(k)
\right) \\
\text{s.t.} \quad
& \eqref{eq:p1-p2} \text{--} \eqref{eq:physical-bat-constraints}, \\
& \eqref{eq:pg-exp-split} \text{--} \eqref{eq:pg-des-split}, \quad \forall k \in \mathcal{K}. \nonumber
\end{aligned}
\end{equation}
By splitting the expected grid power into the two parts $p_G^{exp, buy}(k)$ and $p_G^{exp, sell}(k)$, we enable the usage of different import and export prices $c^{{buy}}(k)$ and $c^{{sell}}(k)$ respectively.
In general, Equations \eqref{eq:pb-exp} and \eqref{eq:pg-exp-split} are not trivial to be evaluated because of the integrals. 
Later in \autoref{sec:model-simulation}, we assume that net-load uncertainty $f_{P_L}$ follows a bimodal Gaussian distribution, which allows us to find closed-form expressions for all integrals. Note however, that this assumption is not a necessity and similar integrals can be approximated without finding a closed-form expression, see our previous work \cite{pinter_probabilistic_2025}.

With that, we obtain a nonlinear, nonconvex stochastic optimization problem that captures uncertainty in both battery and grid power through mixed random variables. Our framework does not output single battery power setpoints, but rather yields battery scheduling policies consisting of a battery power interval (restricted by $\ubar{p}_B$ and $\Bar{p}_B$) and a desired grid power $p_G^{des}$.
In the following, we evaluate our theory by implementing it as a stochastic MPC and simulating its performance across 15 buildings over 5 consecutive months.

\section{Model Simulation}
\label{sec:model-simulation}

To validate our theory, we conduct simulations using real-world data. \autoref{sec:data} introduces the selected dataset, \autoref{sec:forecasts} explains forecast generation, \autoref{sec:prices} specifies the electricity price structure, and \autoref{sec:experiment-configuration} contains additional configuration details and benchmarks. 

The model is implemented in Python using the optimization modeling language Pyomo \cite{bynum_pyomo_2021}. The solution time per optimization iteration is on the order of a few seconds.

\subsection{Data}
\label{sec:data}

Our simulations are based on the dataset presented in \cite{schlemminger_dataset_2022}, which contains high-resolution consumption measurements from residential buildings of a district in northern Germany from 2018 to 2020. We select 15 buildings with the highest data availability; all are equipped with a heat pump, and none has a PV system installed.
To obtain realistic net-load data with PV characteristics, we map PV measurements from three nearby PV sites included in the dataset onto the buildings.
We assign a number of PV panels according to building size, ranging from 16 to 38 panels per building.
The resulting dataset spans two years and five months of real-world net-load data for 15 buildings.

\subsection{Forecasts}
\label{sec:forecasts}

The proposed approach requires probabilistic net-load forecasts. We employ a Kolmogorov–Arnold Network forecasting model \cite{liu_kan_2025} in an autoregressive manner. The model relies solely on historical net-load data without exogenous variables such as weather forecasts or calendar features. The forecasting model generates 24-hour ahead forecasts at hourly resolution. For each hour of a 24-hour forecast, the model outputs 99 point predictions representing the 1\% to 99\% quantiles forming the probabilistic forecast. 

We split the dataset chronologically for time-series cross validation. The first year serves as the training set, the second year as the validation set, and the final five months as the test period. 
For each building, we train an individual local model by evaluating 125 hyperparameter configurations using the Tree-Structured Parzen Estimator algorithm \cite{watanabe_tree-structured_2025}, which samples configurations from a predefined grid and selects the one with the best validation performance. We then generate rolling 24-hour forecasts for the test period.

\begin{figure}[h]
    \centering
    \includegraphics[width=0.93\linewidth]{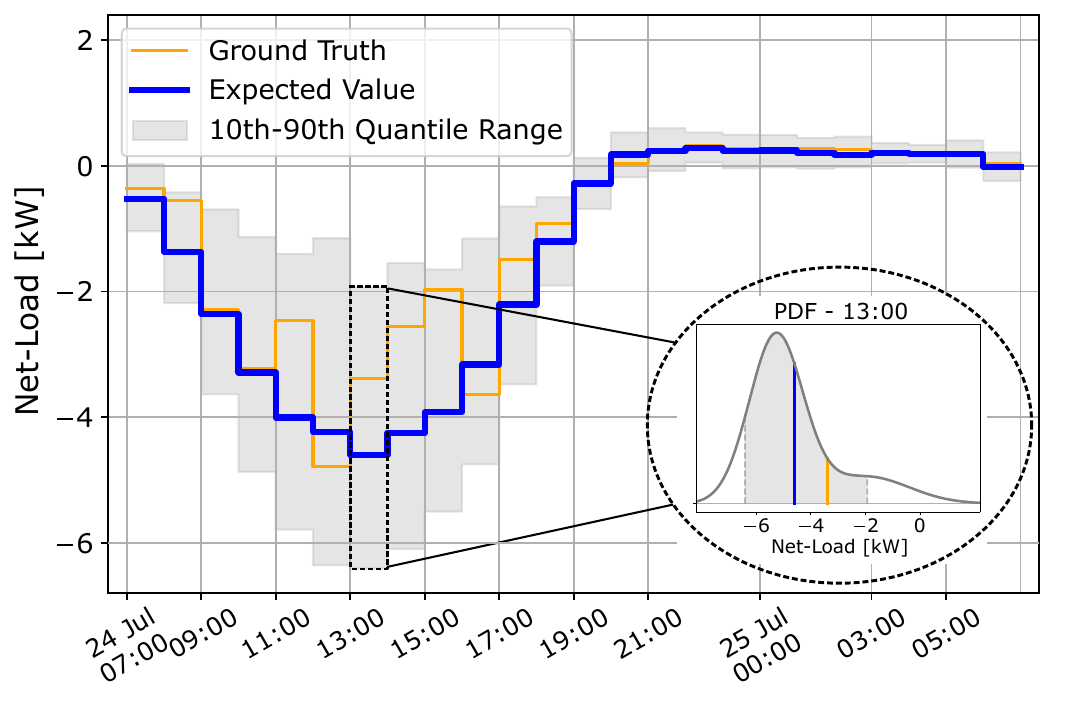}
    \caption{Probabilistic net-load forecast. The lower-right panel visualizes the forecasted PDF at 13:00, including the respective quantiles, the expected value, and the ground truth.}
    \label{fig:prob-forecast}
\end{figure}

To incorporate the probabilistic forecasts into the proposed optimization framework, they must be available in parametric form. We assume that hourly net-load uncertainty can be represented by a continuous random variable whose PDF is the weighted sum of two normal distributions:
\begin{equation}
    f_{P_L}(z) = \frac{1}{\sqrt{2\pi}} \left( \frac{\omega_1}{\sigma_1} e^{\frac{-(z-\mu_1)^2}{2\sigma_1^2}} + \frac{\omega_2}{\sigma_2} e^{\frac{-(z-\mu_2)^2}{2\sigma_2^2}} \right).
\end{equation}
The parameters $\omega_i$, $\mu_i$, and $\sigma_i$ denote the weight, mean and standard deviation of the $i$-th normal component, respectively. They are fitted separately for each hour based on the quantile forecasts. We select this form because it provides a good compromise between computational tractability in the downstream optimization and the ability to represent skewed or otherwise asymmetric net-load behavior.
\autoref{fig:prob-forecast} displays a probabilistic net-load forecast for a typical summer day. 

\subsection{Electricity Price}
\label{sec:prices}

The prices used in this study are derived from EPEX Spot day-ahead prices for the bidding zone Germany-Luxemburg \cite{entsoe}, which are assumed to be fully known in advance. To reflect end-user prices, we convert the wholesale prices of 2025 to retail import prices by adding a constant offset that captures taxes, levies, grid charges and retailer margins. The offset is chosen such that the average price of a year is around 0.4 €/kWh, which corresponds to the average residential electricity price in Germany in 2025 \cite{bundesnetzagentur_preisbestandteile_tarife}.

\begin{figure}[h]
    \centering
    \includegraphics[width=0.8\linewidth]{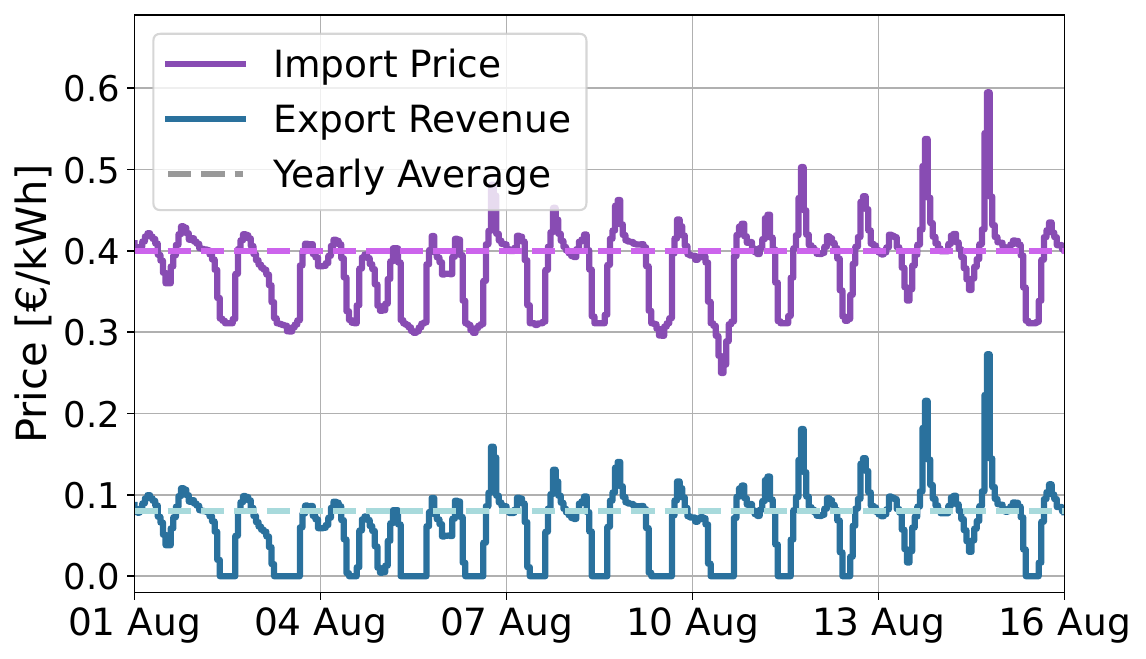}
    \caption{Exemplary two-week period of selected import and export tariffs based on 2025 wholesale electricity prices.}
    \label{fig:prices}
\end{figure}

Export remuneration is obtained analogously. We shift the EPEX wholesale price to an average of 0.08 €/kWh, which roughly matches residential feed-in compensation in Germany as of 2025. Additionally, the export price is clipped at zero to exclude negative feed-in remuneration (i.e., households paying for export). \autoref{fig:prices} illustrates the resulting import costs and export revenues for a two-week period.

\subsection{Experiment Configuration}
\label{sec:experiment-configuration}

To evaluate the proposed approach, battery operation is simulated for 15 residential buildings over five consecutive months of the test set (mid-July to mid-December) using hourly resolution. Battery parameters are selected according to \cite{byd-box}, see \autoref{tab:battery-params}.
The proposed mixed-random-variable-based stochastic MPC method using the Fixed-Grid strategy (\textbf{SMPC-FG}) is compared with five control strategies:
\begin{itemize}
    \item \textbf{MPC with Ideal Forecasts (MPC-Ideal):} MPC optimization based on perfect 24-hour net-load forecasts. 

    \item \textbf{MPC with Fixed Battery (MPC-FB):} MPC optimization based on a deterministic forecast. The first battery setpoint of an MPC iteration is executed, deviations are compensated by the grid. 

    \item \textbf{Stochastic MPC with Fixed Battery (SMPC-FB):} MPC optimization based on a probabilistic forecast oriented at the theory presented in \cite{su_optimal_2021}. The first battery setpoint is executed, deviations are compensated by the grid.

    \item \textbf{MPC with Fixed Grid (MPC-FG):} MPC optimization based on a deterministic forecast. The first grid setpoint is executed, deviations are compensated by the battery (if physically feasible).

    \item \textbf{Rule-Based Control (RBC):} No optimization and no forecast. The battery charges during PV surplus and discharges during load surplus to minimize grid exchange while respecting physical battery constraints. No price information is considered.
\end{itemize}
For deterministic forecasts, we use the expected value of the probabilistic forecasts. By doing so, all forecast-based methods rely on the same underlying forecast information. 

\begin{table}[h]
\centering
\caption{Selected battery specifications.}
\label{tab:battery-params}
\setlength{\tabcolsep}{4.5pt} 
\begin{tabular}{c c c c c}
\toprule
$e^{\min}$ [kWh] & $e^{\max}$ [kWh] & $p_{B}^{\min}$ [kW] & $p_{B}^{\max}$ [kW] & $\eta^{ch/dis}$ [\%] \\
\midrule
0.0 & 7.68 & -5.12 & 5.12 & 98 \\
\bottomrule
\end{tabular}
\end{table}

\section{Results}
\label{sec:results}

The results of the performed simulations are summarized in \autoref{tab:results}. The table contains imported, and exported energy, and their respective costs. Total Costs denotes the sum of import costs and export revenues and is the main performance indicator. Regret describes the percentage increase in Total Costs compared to the MPC-Ideal benchmark. Finally, Rank denotes the average placement of the respective model over the 15 buildings. For example, MPC-Ideal has a Rank of 1.00, which indicates that MPC-Ideal performed better than the other control approaches for all 15 buildings.

\begin{table}[ht]
\centering
\caption{Average Results Across 15 Buildings}
\label{tab:results}
\scriptsize
\setlength{\tabcolsep}{2.5pt}
\begin{tabular}{l c c c c c c c}
\toprule
\textbf{Model} & \textbf{Imp.} & \textbf{Imp. Costs} & \textbf{Exp.} 
               & \textbf{Exp. Rev.}
               & \textbf{Total Costs} & \textbf{Regret} & \textbf{Rank} \\
              & [kWh] & [€] & [kWh] & [€] & [€] & [\%] & [-] \\
\midrule
MPC-Ideal & 1215.16 & 471.31 & 1171.67 & 80.71 & 387.60 & - & 1.00 \\
\cmidrule(lr){1-8}
\textbf{SMPC-FG} & 1273.15 & 496.84 & 1230.96 & 83.89 & \textbf{410.51} & \textbf{6.8} & 2.00 \\
MPC-FG & 1334.36 & 520.36 & 1287.62 & 86.67 & 431.07 & 12.4 & 3.27 \\
RBC & \textbf{1184.56} & \textbf{487.10} & 1176.76 & 38.79 & 448.30 & 21.1 & 4.27 \\
SMPC-FB & 1374.87 & 547.84 & 1341.49 & 94.35 & 451.54 & 19.7 & 4.53 \\
MPC-FB & 1529.59 & 598.65 & \textbf{1488.81} & \textbf{109.35} & 486.66 & 30.6 & 5.93 \\
\bottomrule
\end{tabular}
\begin{tablenotes}
\scriptsize
\item MPC: Model Predictive Control; SMPC: Stochastic MPC; FG: Fixed Grid; FB: Fixed Battery; RBC: Rule-Based Control; SMPC-FG is our contribution. 
\end{tablenotes}
\end{table}

The proposed SMPC-FG approach based on mixed random variables achieves on average a Total Cost of €410.51 with only 6.8\% regret compared to MPC-Ideal and ranks consistently first across all 15 buildings (excluding the MPC based on ideal forecasts). This demonstrates that explicitly modeling the probabilistic grid and battery power interactions through mixed random variables significantly improves performance over deterministic approaches. In comparison, the deterministic MPC-FG achieves 12.4\% regret, nearly double that of SMPC-FG, highlighting the value of modeling uncertainties.

\begin{figure}[h]
    \centering
    \includegraphics[width=0.945\linewidth]{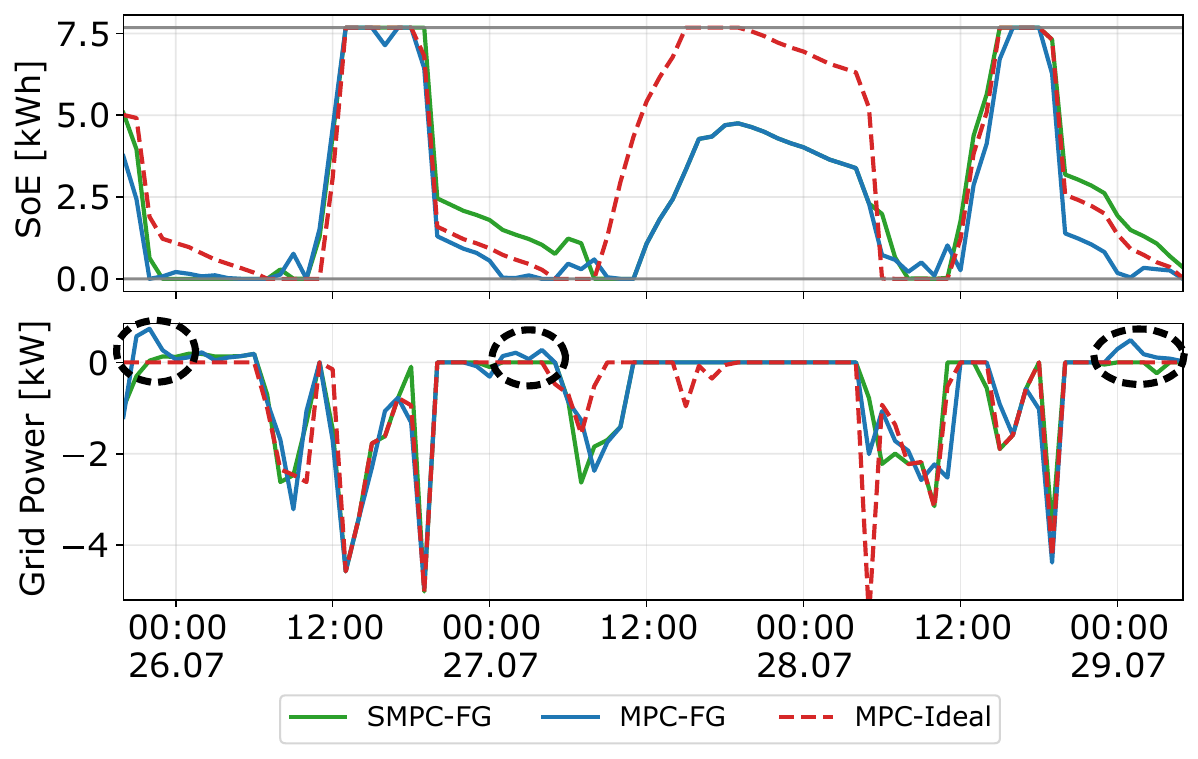}
    \caption{Operation of one building over three summer days. Under the deterministic MPC-FG, the fully discharged battery leads to grid imports at several time steps, causing unnecessary additional costs (black circles).}
    \label{fig:sfh28-operation}
\end{figure}

Furthermore, the performance gap between Fixed-Grid and Fixed-Battery strategies is substantial. MPC-FG outperforms MPC-FB by 18.2 percentage points in regret (12.4\% vs. 30.6\%), and SMPC-FG outperforms SMPC-FB by 12.9 percentage points (6.8\% vs. 19.7\%). 
This is because Fixed-Battery strategies lead to avoidable import-then-export situations, which increase costs under the asymmetric price structure in residential settings. This can also be seen in \autoref{fig:sfh28-operation}, where the behavior of different controllers for one building over three days in July is shown. Overall, our probabilistic SMPC-FG approach and its deterministic counterpart MPC-FG show a similar operating pattern, as both rely on the same forecasts. The main difference is that SMPC-FG behaves more conservatively by discharging the battery less at the beginning of the night. While this reduces export revenues (see \autoref{tab:results}), it preserves flexibility to cover potentially higher demand later in the night. In contrast, MPC-FG depletes the battery more aggressively, which leads to several instances of costly grid imports.



\section{Discussion}
\label{sec:discussion}



Overall, our results emphasize the importance of downstream control in residential energy management, as the Fixed-Grid variants perform substantially better. In many conventional MPC setups, however, this downstream layer is not explicitly considered during optimization. Instead, the MPC computes an optimal decision first, and a subsequent controller is applied afterward as a separate layer. As a result, the MPC decisions do not account for how this downstream controller will respond under different uncertainty realizations. In contrast, our mixed-random-variable-based formulation explicitly incorporates this subsequent control behavior into the optimization. Rather than only providing a target grid or battery setpoint, our approach also captures how different uncertainty realizations should be handled through different charging or discharging actions. This additional information enables more informed decisions as well as more informed downstream control and leads to the improved financial performance observed in our results.

\section{Conclusion}
\label{sec:conclusion}

In this paper, we highlight the benefits of mixed random variables for uncertainty-aware control in residential energy management. Our results show that incorporating mixed random variables into a stochastic MPC framework can improve decision-making under uncertainty and achieve lower electricity costs than other deterministic and probabilistic formulations. A key advantage of the proposed approach is that we do not schedule individual setpoints, but rather battery power intervals, which provide additional flexibility as uncertainty unfolds. These findings highlight the potential of mixed random variables as a promising yet still largely unexplored tool for decision-making under uncertainty as well as uncertainty modeling in general. By demonstrating their use in this work, we hope to encourage future research on their applicability in other problem settings.

\section*{Acknowledgments}
The authors thank the Helmholtz Association for their support under the “Energy System Design” program.

\section*{Declaration}
The authors used Perplexity during writing and coding to improve clarity and fluency. After the usage, the authors reviewed and edited the content as needed and take full responsibility for the content of the published article. 


\bibliographystyle{IEEEtran}
\bibliography{IEEEabrv,references_jp}

\end{document}